\date{January 16, 2001}
\title{percolation on finite graphs}
\author{Itai Benjamini}

\documentclass[12pt,naturalnames]{article}
\usepackage{amsmath}
\usepackage{amsthm}
\usepackage{amsfonts}
\usepackage{graphicx}
\newif\ifhyper\IfFileExists{hyperref.sty}{\hypertrue}{\hyperfalse}
\ifhyper\usepackage{hyperref}\fi

\newif\ifdraft
\drafttrue
\numberwithin{equation}{section}
\numberwithin{figure}{section}

\newtheorem{theorem}{Theorem}
\numberwithin{theorem}{section}

\newtheorem{proposition}[theorem]{Proposition}
\newtheorem{conjecture}[theorem]{Conjecture}

\newtheorem{question}[theorem]{Question}

\newcommand{\Z}{\mathbb{Z}}

\def \P {{\bf P}}

\def \proof {{ \medbreak \noindent {\bf Proof.} }}

\def \_reg {\rightarrow_{\bf reg}}

\def\maxdeg/{\Delta}

\begin{document}
\maketitle

\begin {abstract}
Several questions and few answers regarding percolation on finite graphs are
presented.
\end {abstract}



The following is a note regarding the asymptotic study of percolation on finite
transitive graphs.
On the one hand, the theory of percolation on infinite graphs is rather developed,
although still with many open problems (See \cite{L}). On the other hand random graphs
were deeply studied (see \cite{JLR}). Finite transitive graphs are somewhat in the
middle.
Though a big part of the picture is similar to the percolation on infinite graphs
theory, some qualitative phenomena are different and provides new challenges.
\medskip

Recall that a graph $G$ is transitive iff for every two vertices $v,u \in G$,
there is an isomorphism of $G$, mapping $v$ to $u$.
See \cite{BS} or \cite{L} for further background and definitions.
\medskip

Given a finite graph $G$, delete edges from $G$ independently, each with probability
$1-p$. Denote by $p^G= p^G_{1/2}$, the threshold probability for having a connected
component of size $|G|/2$. I.e.,
$$
\P_{p^G}\big( \mbox{There is a connected component of size} \geq |G|/2 \big) = 1/2.
$$

Note that for transitive graphs the threshold  for having a large connected component
is sharp \cite{FK}. Write $p^G_{\alpha}$ for the threshold probability for having a
component of size $\alpha|G|$, with probability $1/2$.

\section{Giant component}

A giant component theorem will be proven under a weak concentration of measure
assumption.
A family of graphs $G_n$, $|G_n|\nearrow \infty$,
admits a weak concentration of measure, if for every $C > 0$, and every
two sets of vertices $A_n, B_n \subset G_n, |A_n|, |B_n| \geq C|G_n|$,

$$
d(A_n, B_n) = o(\mbox{diameter} (G_n)).
$$

Assume also that $\mbox{diameter} (G_n) \leq C \log |G_n|$ and that $K$ is a
uniform upper bound on the degrees of $G_n$.

\begin{theorem}
Let $G_n$ is a sequence of graphs satisfying the conditions above,
then for every $\epsilon > 0$ and $0<C<1$,
$$
\P_{p^{G_n} + \epsilon}\big( \mbox{There is a unique connected component of size} \geq
C|G| \big) \longrightarrow 1,
$$
as $n$ goes to infinity.
\end{theorem}


\proof
It is enough to show that given any two set
$A_n, B_n \subset G_n, |A_n|, |B_n| \geq C|G_n|$,  $C>0$, and $\epsilon >0$,
if we do $\epsilon$ percolation on $G_n$, there will be an open path from
$A_n$ to $B_n$,with probability going to $1$, with $n$.
We will show that
there are an exponential number, in the $\mbox{diameter} (G)$, edge disjoint paths
of length
$o(\mbox{diameter} (G))$ connecting $A_n$ to $B_n$. Indeed for every vertex in $A_n$
look at it's shortest paths
to $B_n$, as $n$ grows for most of the vertices of $A_n$ their shortest path to
$B_n$, has length which is $o(\mbox{diameter} (G))$, so we get a set of paths of
size bigger than $|A_n|/2$, with length much smaller than the diameter, denoted
by $l_n$.
Now for a fixed $l$, the number of paths of length $l$, that intersects a paths of
length $m$ in $G_n$, is bounded by $m \times K^l$.
Which implies that there are more than,
$$
  {|A_n| \over 2l_n \times K^{l_n}}
$$
disjoint paths of length smaller than $l_n$ between $A_n$ and $B_n$.
Each path has probability $ \geq \epsilon^{l_n}$ to be open, independently
from the other paths. The probability of having no open paths is bounded then by,
$$
(1- \epsilon^{l_n})^{{|A_n| \over 2l_n \times K^{l_n}}}.
$$
Pick $n$ large to finish.
\qed
\medskip

The above concentration condition is clearly satisfied by expanders and other
families of graphs, but not by all transitive graphs of logarithmic diameter.
E.g., let $G_n'$ be the product  $G_n \times C_{\mbox{diameter} (G_n)}$,
where $C_n$ is a cycle of length $n$, and $G_n$ satisfy the above assumptions.






\medskip
One would expect that concentration of measure is not necessary,
but then we assume that the graphs are transitive.

\begin{conjecture}
Assume $G_n$ is a sequence of finite transitive graphs with $|G_n|\nearrow \infty$.
Then for every $\epsilon > 0$ and $C > 0$,
$$
\P_{p^{G_n} + \epsilon}\big( \mbox{There is a unique connected component of size} \geq
C|G| \big) \longrightarrow 1,
$$
as $n$ goes to infinity.
\end{conjecture}

The conjecture is true for a large Euclidean torus by the Grimmett-Marstrand theorem
\cite{GM}.

The following proposition (due to Oded) is of independent interest.
Consider percolation on a finite Cayley graph,
\begin{proposition}

Assume
$$
\P_p ( v \mbox{  is in a connected component of size } \geq C|G| ) = q.
$$
Then for every $ u \in G$,
$$
\P_p ( v \mbox{ and } u \mbox{ are in the same connected component})
\geq (Cq/2)^{6/qC }.
$$
\end{proposition}

\proof
The expected number of vertices $v$ is connected to, is bigger or equal to $ qC|G|$.
Thus there
are more than $qC|G|/2$ vertices $u$ which are connected to $v$ with probability
$\geq qC/2$. Let $S$ be the set of these vertices. $S$ is symmetric and generating,
since it contains $v$ neighbors.
By FKG it is enough to show that the diameter of the group using $S$ as a generating
set is bounded by $ 6/qC$.
Assume $a_1a_2.....a_l$ is a shortest word, in the generators $S$, representing an
element of the group. It follows that for any $k$ and
$l$, divisible by $3$, $a_1a_2...a_kg \neq a_1a_2...a_lg'$ for any $g, g' \in S$.
Otherwise one can reduced the length of the word. We get that
$l/3 \times |S| \leq |G|$ i.e. $l \leq 6/qC$.
\qed

\section{Threshold bounds}

In light of the theorem and the conjecture,
(which suggests a further study in the spirit of random
graphs theory, e.g., components size and scaling around the birth of the giant
component etc. see \cite{JLR}), it will be useful to have
upper bounds on $p^G$, guaranteeing the threshold is not converging to $1$.
The techniques mentioned in this section are similar to the infinite
graphs case.

\cite{BB} gives a bound in terms of the length of the longest relation in the
presentation of the group. Theorem 2 in \cite{BS} can be adapted to finite graphs,
and will provide a bound in terms of the isoperimetric
constant of the graph, guaranteeing that the threshold is not converging to $1$
for expanders. These bounds together do not provide a sharp condition
for boundedness away from $1$.
We suspect the following condition is sharp.

\begin{conjecture}
\label{trash}
Assume that $G$ is a transitive graph.
There is a constant $C < 1$, so that if
$$
\mbox{diameter} (G) < {|G| \over \log|G|},
$$
then $p^G < C$.
\end{conjecture}

Maybe the classification theorem for simple groups can be useful in attacking the
conjecture, even just for Cayley graphs. A good start might be Cayley graphs of $S_n$.

One way to prove an upper bound on the threshold for a graph $G$, is to find a
known percolating graph
spanning a linear portion of $G$. Another is to use measures on paths with
exponential intersection tail ($EIT$).

Say a  graph $G$, has the $EIT(c), 0 <c < 1$, property iff
for any two vertices $u,v \in G$,  there is a set of paths, $S_{uv}$, from $v$
to $u$, and a measure on the paths so that if you pick two paths independently at
random from $S_{uv}$, then the probability that the two paths will have more than $k$
intersections decays faster than $c^k$.

A second moment argument (done in the infinite graph setting  in \cite{BPP})
implies that if $G$ has the $EIT(c)$
property, for some $c<1$, the threshold for giant component will be bounded away
from $1$.

Oriented paths on $\Z^4$, provide a family of such paths for $(\Z/n\Z)^4$.
( For $(\Z/n\Z)^3$, it is  harder to construct paths with
exponential intersection tail, and appears in \cite{BPP}.

\begin{conjecture}
\label{eit}
Let $G$ be a finite transitive graph. There is $0 < c <1$, such that if

$$
\mbox{diameter} (G) < |G|^{1/3},
$$
then $G$ admits the $EIT(c)$ property.
\end{conjecture}

Adaptation of the argument in \cite{L1} gives that if
$\mbox{diameter} (G) < C \log |G|$, then
$G$ admits the $EIT(c)$ property (with $c$ depending on $C$).

If conjecture~\ref{eit} is true, it will also imply the following a bit
unrelated conjecture.

\begin{conjecture}
\label{electric}
Let $G$ be a finite transitive graph. There is $ R < \infty$, such that if

$$
\mbox{diameter} (G) < |G|^{1/3},
$$
then  the effective electric resistance between any two vertices in $G$ is smaller
than $R$.
\end{conjecture}



A natural conjecture sharpening theorem 2 of \cite{BS} is
as follows. Given a graph $G$, let
$$
\delta(G) = \sup\{\alpha | |\partial S| \geq |S|^{\alpha}, \mbox{for all } S \subset G,
|S| \leq |G|/2  \}.
$$
\begin{conjecture}
If $\delta(G) > \delta > 0$, then $p^G < f(\delta) < 1$.
\end{conjecture}

\section{Percolation implies giant component?}

Conjecture~\ref{trash} above is a finite analog of the still open conjecture
from \cite{BS}, that the critical
percolation probability is strictly smaller than $1$ on infinite Cayley graphs
of groups which are not finite extension of $\Z$. We now move to discuss the
percolation probability function.

Note that in \cite{BB} it was conjectured that there are no infinite connected
component for critical percolation on infinite transitive graphs and one expect
that the percolation probability function should be continuous, (see \cite{S}).
It is therefore natural to conjecture that for a sequence of finite transitive
graphs $G_n \longrightarrow G$, then for any $ 0 < \alpha_1 < \alpha_2 <1$,
$$
p^{G_n}_{\alpha_1} - p^{G_n}_{\alpha_2} \rightarrow c >0.
$$
I.e. the thresholds for connected component of different proportion stay apart.
In the last section we discussed bounds on $p^G$, i.e. birth of a component of half
the size. In this section we discuss the birth of components of linear size.

Assume the transitive graphs $G_n$ converges weakly to the infinite graph $G$
and $\{p^{G_n}\}_n$ are bounded away from $1$. Let
$\theta_n(p)$ be the probability a fixed vertex is in the largest component, and
$\theta(p)$ be the probability a fixed vertex in $G$ is in an infinite cluster.
As usual $p_c(G)$ is critical probability for percolation on $G$.

\begin{question}
\label{percgiant}
Does
$$
\theta_n(p) \longrightarrow \theta(p)?
$$
And in particular
$$
\lim_{\alpha \rightarrow 0} {\lim_n p^{G_n}_{\alpha}} = p_c(G)?
$$
\end{question}

For example, assume $G_n$ has growing girth and all vertices has the same degree,
say $3$, and $\{p^{G_n}\}_n$ are bounded away from $1$.
Hence $G_n \longrightarrow T_3$, the $3$-regular tree.
Recall that $ p_c(T_3) = 1/2$, that is the critical
probability for percolation on $T_3$ is $1/2$. We have,
\begin{conjecture}
For every $ p > 1/2$ there is  some $C_p >0$, such that
$$
\P_p \big( \mbox{There is a unique connected component of size } C_p|G_n| \big) \longrightarrow
1,
$$
as $ n \longrightarrow \infty$.
\end{conjecture}

An example that supports the conjecture is as follows, start with a cycle of $n$
vertices. Pick uniformly a random matching between the vertices on the cycle. Add
an edge between any matched pair. Around a typical vertex this graph (a typical sample)
looks like a regular tree for more than $1/10 \log n$ generations, and indeed as $n$
grows the threshold for the existence of a giant component goes to $1/2$.

Of independent interest, we
expect the largest components at the threshold to be of size $O(n^{2/3})$ as in the
critical random graphs $G(n, 1/n)$ and probably  in critical percolation on
expanders. In particular we conjecture that for percolation with $p =(k-1)^{-1}$ on a
random $k$-regular graph (that's the threshold). The largest components is of size
$O(n^{2/3})$.

Another supporting example is the theorem in \cite{AKS} asserting
that the threshold for a giant component on the hypercube $\{0,1\}^d$, is $1/d$.
Conjecture~\ref{percgiant}
is true for a large Euclidean torus by the Grimmett-Marstrand theorem
\cite{GM}.

Note that, if $B_n$ are growing balls in a nonamenable Cayley graph $G$, assuming
the conjecture (from \cite{BS}) that $p_c(G) < p_u(G)$, then the threshold for
connected component of linear size in $B_n$ will converge to $p_u$ rather than $p_c$.

Intuitively the content of the last two conjectures is as follows:
On infinite transitive graph percolation means having unbounded connected components,
and uniqueness of the infinite cluster is the analog of existence of  a giant
component. We doubt,
whether on finite transitive graphs component of large diameter can occur and moreover
contain portion of the vertices, without
having a giant component. One informal reason we have, is that unlike infinite graphs
(e.g. regular trees) finite transitive graphs do not have bottlenecks, unless they
contain an irreducible cycle, proportional in size to the graphs, (as in $C_n \times G$,
where $C_n$ are growing cycles and $G$ is a fixed finite graph). An argument using
flows (see for instance \cite{SC}) gives that the Cheeger constant of a finite
transitive graph, $h(G)$, satisfies
$h(G) \geq \mbox{diameter}^{-1}(G)$.

\medskip

\noindent
{\bf Remark:} (On limits of transitive graphs, by Russ and Yuval)
Every weak limit of finite transitive graphs (or even of unimodular
infinite transitive graphs) is unimodular transitive. This is from the
fact that there is a unimodular group acting transitively on the limit
graph, namely, the group of automorphisms that arise as limits of
automorphisms. This group is unimodular by the criterion of Trofimov \cite{T}.

Let us end by saying that there are further directions to pursue, e.g., compare
mixing times of random walks on the giant component to mixing times on the graph.
Maybe one should start with a further  study of the geometry of finite
transitive graphs e.g., what asymptotic diameters possible for finite
transitive graphs, or cutsets structure.

\noindent{\bf Acknowledgments:} Thanks to Russell Lyons, Igor Pak, Yuval Peres and Oded Schramm
for useful discussions.

\medskip
\noindent


\begin{thebibliography}{AAA}

\bibitem{AKS}
Ajtai, M., Komlos, J. and Szemerédi, E.,
Largest random component of a $k$-cube. Combinatorica 2 (1982), no. 1, 1--7.

\bibitem{BB}
Babson, E. and Benjamini, I., Cut sets and normed cohomology with applications to percolation.
Proc. Amer. Math. Soc. 127 (1999), no. 2, 589--597.


\bibitem{BPP}
Benjamini, I., Pemantle, R. and Peres, Y., Unpredictable paths and percolation.
Ann. Probab. 26 (1998), no. 3, 1198--1211.

\bibitem{BS}
Benjamini, I. and Schramm, O., Percolation beyond $\bold Z\sp d$, many questions and a few answers.
Electron. Comm. Probab. 1 (1996), no. 8, 71--82 (electronic).


\bibitem{FK}
Friedgut, E. and Kalai, G., Every monotone graph property has a sharp threshold.
Proc. Amer. Math. Soc. 124 (1996), no. 10, 2993--3002.


\bibitem{GM}
Grimmett, G. R. and Marstrand, J. M., The supercritical phase of percolation is well
behaved.
Proc. Roy. Soc. London Ser. A 430 (1990), no. 1879, 439--457.

\bibitem{JLR}
Janson, S., \L uczak, T. and Rucinski, A., Random graphs. Wiley-Interscience Series in Discrete Mathematics and Optimization.
Wiley-Interscience, New York, 2000. xii+333 pp.

\bibitem{L1}
Lyons, R., Random walks and the growth of groups.
C. R. Acad. Sci. Paris Sér. I Math. 320 (1995), no. 11, 1361--1366.

\bibitem{L}
Lyons, R., Phase transitions on nonamenable graphs.
Probabilistic techniques in equilibrium and nonequilibrium statistical physics.
J. Math. Phys. 41 (2000), no. 3, 1099--1126.


\bibitem{SC}
Saloff-Coste, L., Lectures on finite Markov chains. Lectures on probability theory and statistics
(Saint-Flour, 1996), 301--413, Lecture Notes in Math., 1665, Springer, Berlin, 1997.

\bibitem{S}
Schonmann, R. H., Stability of infinite clusters in supercritical percolation.
Probab. Theory Related Fields 113 (1999), no. 2, 287--300.

\bibitem{T}
Trofimov, V. I., Groups of automorphisms of graphs as topological groups.
(Russian) Mat. Zametki 38 (1985), no. 3, 378--385, 476.

\end{thebibliography}
\end{document}